\documentclass{amsproc}

\usepackage{url}

\newtheorem{theorem}{Theorem}[section]

\newtheorem{proposition}[theorem]{Proposition}
\theoremstyle{definition}
\newtheorem{definition}[theorem]{Definition}

\theoremstyle{remark}

\numberwithin{equation}{section}

\newcommand{\Irr}[1]{\operatorname{Irr}(#1)}

\begin{document}

\title[  Geck-Rouquier classification of simple modules for Hecke algebras    ]{JMMO Fock space and Geck-Rouquier classification of simple modules for Hecke algebras}

\author{Nicolas Jacon}

\address{Laboratoire de Math\'ematiques Nicolas Oresme, Universit\'e de Caen, BP 5186, F 14032 Caen Cedex, France.}

\email{jacon@math.unicaen.fr}



\begin{abstract}
Using Lusztig $a$-function, M.Geck and R.Rouquier have recently proved the existence of a canonical set $\mathcal{B}$ in natural bijection with the set of simple modules for Hecke algebras. In this paper, we recall the definition of this set and we report recent results which show that the definition of canonical basic set can be extended to the case of Ariki-Koike algebras. Moreover, we  give an explicit description of this set for all Hecke algebras and for all Ariki-Koike algebras.
\end{abstract}

\maketitle

\section{Introduction}
Let $W$ be a finite Weyl group with set of simple reflections $S\subset{W}$, let $v$ be an indeterminate, $u=v^2$ and let $A=\mathbb{C}[v,v^{-1}]$. The Hecke algebra $H$ of $W$ over $A$ is  the associative $A$-algebra with  basis $\{T_w\ |\ w\in{W}\}$ and the multiplication between two elements of the basis is determined by the following rule. Let $s\in{S}$ and $w\in{W}$, then:
$$ T_s T_w=
    \begin{cases}
        T_{sw} & \textrm{if } l(sw)>l(w), \\
       uT_{sw}+(u-1)T_w & \textrm{if } l(sw)<l(w),\\\end{cases}$$
where $l$ is the usual length function of $W$.  Such algebras play an important role, for example in the representation theory of finite groups of Lie type (see \cite{DG} and \cite{Gb}) or in the theory of knots and links (see \cite[Chapter 4]{GP}). 

Let $K$ be the field of fractions of $A$ and let $\theta : A \to \mathbb{C}$ be a homomorphism into the field of complex numbers. Let $H_K:=K\otimes_A H$ and let $H_{\mathbb{C}}:=\mathbb{C} \otimes_A H$. On the one hand, the representation theory of $H_K$ is relatively well-understood: it is known that this is a split semisimple algebra which is isomorphic to the group algebra  $\mathbb{C}[W]$ and its simple modules are in natural bijection with the simple  $\mathbb{C}[W]$-modules. On the other hand, the simple modules of $H_{\mathbb{C}}$ are much more complicated to describe because $H_{\mathbb{C}}$ is not semisimple in general. In fact, this problem is linked to the problem of determining  a map   $d_{\theta}$ between the Grothendieck group $R_0 (H_K)$ of finitely generated $H_K$-modules and the Grothendieck group $R_0 (H_{\mathbb{C}})$ of finitely generated $H_{\mathbb{C}}$-modules. This map  is called the  ``decomposition map'' and it relates the simple $H_K$-modules with the simple $H_{\mathbb{C}}$-modules via a process of modular reduction: $$d_{\theta}:R_0(H_K)\to R_0(H_{\mathbb{C}}).$$
The decomposition maps associated to Hecke algebras of exceptional types are almost all explicitely known, see \cite{Mu}, \cite{Ge1}, \cite{Ge2} and \cite{GL} (in fact, for $W=E_8$, we only have an ``approximation'' of this map, see \cite{Mu}). In this paper, we are mainly interesting about the classical types that is type $A_{n-1}$, type $B_n$ and type $D_n$. 

 For these types, works of Ariki \cite{Ab}, Dipper-James \cite{DJ}, Dipper-James-Murphy \cite{DJM},  and Hu \cite{H} provide a description of the simple $H_{\mathbb{C}}$-modules. In fact, Hecke algebras of type $A_{n-1}$ and $B_n$ are particular cases of Ariki-Koike algebras (or cyclotomic Hecke algebras of type $G(d,1,n)$, see \cite{AK}). In \cite{A2},  Ariki has shown that the computation of the decomposition map for these algebras can be easily deduced from the computation of the Kashiwara-Lusztig  canonical basis of irreducible highest weight $\mathcal{U}_q (\widehat{sl}_e)$-modules. In particular, these results lead to a parametrization of the simple modules for Ariki-Koike algebras (and so, for Hecke algebras of type $A_{n-1}$ and $B_n$) using a class of  multipartitions which appears in the crystal graph theory of $\mathcal{U}_q (\widehat{sl}_e)$-modules, namely the Kleshchev multipartitions. For type $D_n$, Hu has obtained a classification of the simple modules by using the fact that the Hecke algebra of type $D_n$ can be seen as a subalgebra of a Hecke algebra of type $B_n$ (with unequal parameters). One of the main problem is that, for type $B_n$ and $D_n$, these results lead to a recursive parametrization of the simple modules.

In \cite{Gk} and in \cite{GR}, Geck and Rouquier have given another approch for the description of the simple $H_{\mathbb{C}}$-modules. They have    shown the existence  of  a canonical set  $\mathcal{B}\subset{\Irr{H_K}}$ by using Lusztig  $a$-fonction and Kazhdan-Lusztig theory. This set is called the ``canonical basic set'' and it is in natural bijection with $\Irr{H_{\mathbb{C}}}$. Hence, it gives a way to classify the simple $H_{\mathbb{C}}$-modules. Moreover, the existence of the canonical basic set implies that the matrix associated to $d_{\theta}$ has a lower triangular shape with $1$ along the diagonal.

The description of the canonical basic set is now complete for all finite Weyl groups $W$ and for all specializations $\theta$. In this paper, we report these recent results which show that this  set can be indexed by another class of multipartitions which also appears in the crystal graph theory of the  $\mathcal{U}_q (\widehat{sl}_e)$-modules. The proof  also requires Ariki's theorem but what we obtain here is a non recursive parametrization of the simple  $H_{\mathbb{C}}$-modules. Moreover, we will see that all these results can be extended to the case of Ariki-Koike algebras even if we don't have Kazhdan-Lusztig type basis for Ariki-Koike algebras.

\section{Representations of semisimple Hecke algebras}

Let $H$ be an Iwahori-Hecke algebra of a finite Weyl group $W$ over $A:=\mathbb{C}[v,v^{-1}]$ as it is defined in the introduction. Let $K=\mathbb{C}(v)$ and let $H_K$ be the corresponding Hecke algebra. Then $A$ is integrally closed in $K$ and  $H_K$ is a split semisimple  algebra. By Tits deformation theorem (see \cite[Theorem 8.1.7]{GP}), $H_K$ is isomorphic to the group algebra $\mathbb{C}[W]$. Hence, the simple  $H_K$-modules are in natural bijection with the simple $\mathbb{C}[W]$-modules. In fact,  for type $A_{n-1}$ and type $B_n$, the simple $H_K$-modules can be explicitly described by using the theory of cellular algebras (see \cite{GrL}) while for type $D_n$, the simple $H_K$-modules are obtained by using Clifford theory. We obtain the following parametrizations for the classical types of Weyl groups:

\subsection{Type $A_{n-1}$} Assume that $W$ is a Weyl group of type $A_{n-1}$. 
\\
\\
\begin{center}
\begin{picture}(200,20)
\put( 40,10){\circle*{5}}
\put( 37,18){$s_1$}
\put( 40,10){\line(1,0){40}}
\put(80,10){\circle*{5}}
\put(77,18){$s_2$}
\put(80,10){\line(1,0){20}}
\put(110,10){\circle{1}}
\put(120,10){\circle{1}}
\put(130,10){\circle{1}}
\put(140,10){\line(1,0){20}}
\put(160,10){\circle*{5}}
\put(157,18){$s_{n-1}$}
\end{picture}
\end{center}
\vspace{0,5cm}
Let $\lambda$ be a partition of rank $n$, then, we can construct a $H$-module $S^{\lambda}$, free over $A$  which is called  a Specht module (see the construction of ``dual Specht modules'' in \cite[Chapter 13]{Ab} in a more general setting). Moreover, we have:
$$\Irr{H_K}=\{S^{\lambda}_K:=K\otimes_A S^{\lambda}\ |\ \lambda\in{\Pi^1_n}\},$$
where we denote by $\Pi^1_n$ the set of partitions of rank $n$.

\subsection{Type $B_n$} Assume that $W$ is a Weyl group of type $B_n$. 
\\
\\
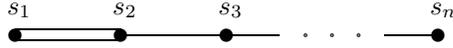
\begin{center}
\begin{picture}(240,20)
\put( 50,10){\circle*{5}}
\put( 47,18){$s_1$}
\put( 50,8){\line(1,0){40}}
\put( 50,12){\line(1,0){40}}
\put( 90,10){\circle*{5}}
\put( 87,18){$s_2$}
\put( 90,10){\line(1,0){40}}
\put(130,10){\circle*{5}}
\put(127,18){$s_3$}
\put(130,10){\line(1,0){20}}
\put(160,10){\circle{1}}
\put(170,10){\circle{1}}
\put(180,10){\circle{1}}
\put(190,10){\line(1,0){20}}
\put(210,10){\circle*{5}}
\put(207,18){$s_{n}$}
\end{picture}
\end{center}
\vspace{0,5cm}
Let $(\lambda^{(0)},\lambda^{(1)})$ be a bi-partition of rank $n$, then, we can construct a $H$-module $S^{(\lambda^{(0)},\lambda^{(1)})}$, free over $A$  which is called  a Specht module. Moreover, we have:
$$\Irr{H_K}=\{S^{(\lambda^{(0)},\lambda^{(1)})}_K:=K\otimes_A S^{(\lambda^{(0)},\lambda^{(1)})} \ |\ (\lambda^{(0)},\lambda^{(1)})\in{\Pi^2_n}\},$$
where we denote by $\Pi^2_n$ the set of bi-partitions of rank $n$.
\subsection{Type $D_n$}\label{dd} Assume that $W$ is a Weyl group of type $D_n$.\\
\\
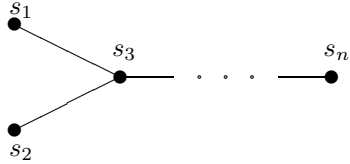
\begin{center}
\begin{picture}(200,20)
\put( 40,30){\circle*{5}}
\put( 38,34){$s_1$}

\put( 40,-10){\circle*{5}}
\put( 38,-20){$s_2$}

\put(80,10){\line(-2,1){40}}
\put(80,10){\line(-2,-1){40}}

\put(80,10){\circle*{5}}
\put(77,18){$s_3$}
\put(80,10){\line(1,0){20}}
\put(110,10){\circle{1}}
\put(120,10){\circle{1}}
\put(130,10){\circle{1}}
\put(140,10){\line(1,0){20}}
\put(160,10){\circle*{5}}
\put(157,18){$s_{n}$}
\end{picture}
\end{center}
\vspace{0,8cm}
Then, $H$ can be seen as a subalgebra of a Hecke algebra $H_1$ of type $B_n$ with unequal parameters (see \cite{Gex} for more details).

 Similary to the equal parameter case, for all $(\lambda^{(0)},\lambda^{(1)})\in{\Pi^2_n}$,   we can construct a $H_1$-module $S^{(\lambda^{(0)},\lambda^{(1)})}$, free over $A$  which is called  a Specht module. We have:
$$\Irr{H_{1,K}}=\{S^{(\lambda^{(0)},\lambda^{(1)})}_K:=K\otimes_A S^{(\lambda^{(0)},\lambda^{(1)})}\ |\ (\lambda^{(0)},\lambda^{(1)})\in{\Pi^2_n}\}.$$
We have an operation of restriction $\operatorname{Res}$ between the set of $H_{1,K}$-modules and the set of $H_K$-modules, for $(\lambda^{(0)},\lambda^{(1)})\in{\Pi_n^2}$:
\begin{itemize}
\item if $\lambda^{(0)}\neq \lambda^{(1)}$, we have $\operatorname{Res}(S^{(\lambda^{(0)},\lambda^{(1)})}_K )\simeq \operatorname{Res}(S^{(\lambda^{(1)},\lambda^{(0)})}_K )$ and the $H_K$-module $V^{[\lambda^{(0)},\lambda^{(1)}]}:=\operatorname{Res}(S^{(\lambda^{(0)},\lambda^{(1)})}_K )$ is a simple $H_K$-module.
\item if $\lambda^{(0)}=\lambda^{(1)}$, we have $\operatorname{Res}(S^{(\lambda^{(0)},\lambda^{(1)})}_K )=V^{[\lambda^{(0)},+]}\oplus V^{[\lambda^{(0)},-]}$ where $V^{[\lambda^{(0)},+]}$ and $V^{[\lambda^{(0)},-]}$ are non isomorphic simple $H_K$-modules.
\end{itemize}
Moreover, we have:

$$\Irr{H_K}=\left\{V^{[\lambda,\mu]}\ |\ \lambda\neq{\mu},\ (\lambda,\mu)\in{\Pi^2_n}\right\}
\bigcup{\left\{V^{[\lambda,\pm]}\ |\ \lambda\in{\Pi^1_{\frac{n}{2}}} \right\}}.$$
\vspace{0,3cm}

Now, we turn to the problem of determining a classification of the set of simple modules for Hecke algebras in the case where $v$ is no longer an indeterminate but a complex number.

\section{Modular representations and canonical basic sets for Hecke algebras}

Let  $\theta : A \to \mathbb{C}$ be a ring homomorphism. We put:
$$\mathcal{O}:=\left\{\frac{f}{g}\ |\ f,g\in{\mathbb{C}[v]},\ g(\theta(v))\neq 0\right\}.$$
 $\mathcal{O}$ is a discrete valuation ring and we have  $A\subset \mathcal{O}$. By \cite[Theorem  7.4.3]{GP}, we obtain a well-defined decomposition map
$$d_{\theta}:R_0(H_K)\to R_0(H_{\mathbb{C}}).$$
where $R_0 (H_K)$ (resp. $R_0 (H_{\mathbb{C}})$) is the Grothendieck group of finitely generated $H_K$-modules (resp. $H_{\mathbb{C}}$-modules). 
This is defined as follows: let $V$ be a simple $H_K$-module. Then, by \cite[\S 7.4]{GP}, there exists a $H_{\mathcal{O}}$-module $\widehat{V}$ such that $K\otimes_{\mathcal{O}}\widehat{V}=V$. By reducing $\widehat V$ modulo the maximal ideal of $\mathcal{O}$, $\mathfrak{m}:=(v-\theta (v))\mathcal{O}$, we obtain a $H_{\mathbb{C}}$-module $\mathbb{C}\otimes_{\mathcal{O}}\widehat V$ .  Then, we put:$$d_{\theta}([V])=[\mathbb{C}\otimes_{\mathcal{O}}\widehat V].$$
$d_{\theta}$ is well-defined and  for $V\in{\Irr{H_K}}$, there exist numbers $(d_{V,M})_{M\in{\Irr{H_{\mathbb{C}}}}}$ such that:
$$d_{\theta}([V])=\sum_{M\in{\Irr{H_{\mathbb{C}}}}} d_{V,M}[M].$$
The matrix $(d_{V,M})_{{V\in{\Irr{H_K}}}\atop{M\in{\Irr{H_{\mathbb{C}}}}}}$ is called the decomposition matrix. For more details about the construction of decomposition maps, even in a more general setting, see \cite{Gb}. 

Now, we will recall results of  Geck and  Rouquier which show that the decomposition map has always a unitriangular shape with one along the diagonal.

Let $\{C_w\}_{w\in{W}}$ be the Kazhdan-Lusztig basis of $H$. For $x,y\in{W}$, the multiplication between two elements of this basis is given by:
$$C_x C_y=\sum_{z\in{W}} h_{x,y,z} C_z$$
where $h_{x,y,z}\in{A}$ for all $z\in{W}$. For any $z\in{W}$, there is a well-defined integer $a(z)\geq 0$ such that 
\begin{align*}
& v^{a(z)}h_{x,y,z}\in{\mathbb{Z}[v]}  \textrm{ for all }x,y\in{W},\\
& v^{a(z)-1}h_{x,y,z}\notin{\mathbb{Z}[v]}  \textrm{ for some }x,y\in{W}.
\end{align*}
We obtain a function which is called the Lusztig  $a$-function:
$$\begin{array}{cccc}
   a : & W & \to & \mathbb{N}\\
     & z & \mapsto & a(z)
\end{array}
$$
Now, following \cite[Lemma 1.9]{L1}, to any $M\in{\Irr{H_{\mathbb{C}}}}$, we can attach an $a$-value $a(M)$ by the requirement that:
\begin{align*}
C_w.M=0 & \textrm{ for all }w\in{W}\textrm{ with } a(w)>a(M),\\
C_w.M\neq 0 & \textrm{ for some }w\in{W}\textrm{ with } a(w)=a(M).
\end{align*}
We can also attach an $a$-value $a(V)$ to any $V\in{\Irr{H_K}}$, in an analogous way. Note that there is  an equivalent definition of the $a$-value of a simple $H_K$-module using the fact that Hecke algebras are symmetric algebras, this will be important in the context of Ariki-Koike algebras where we don't have Kazhdan-Lusztig theory. Let $\tau: H_K\to A$ be the symmetrizing trace of $H$ (see \cite[\S 7.1]{GP}) which is defined by $\tau (T_w)=0$ if $w\neq 1$ and $\tau (T_1)=1$. Then, for each $V\in{\Irr{H_K}}$, there exists a Laurent polynomial $s_V\in A$ such that:
$$\tau=\sum_{V\in{\Irr{H_K}}} \frac{1}{s_V} \chi_V,$$
where $\chi_V$ is the character  of $V\in{\Irr{H_K}}$. $s_V$ is called the Schur element of the simple $H_K$-module $V$. Then, Lusztig (\cite{L2}) has shown that we have:
$$a(V)=\frac{1}{2}\textrm{min}\{l\in{\mathbb{Z}}\ |\ v^ls_V\in{\mathbb{Z}[v]}\}.$$
 We can now give the theorem of existence of the canonical basic set. The main tool of the proof  is the Lusztig  asymptotic algebra.
\begin{theorem}[Geck \cite{Gk}, Geck-R.Rouquier \cite{GR}]\label{basicset} We define the following subset of $\Irr{H_K}$:
$$\mathcal{B}:=\{V\in{\Irr{H_K}}\ |\ d_{V,M}\neq 0\ \textrm{and}\ a(V)=a(M)\ \textrm{for some }M\in{\Irr{H_{\mathbb{C}}}}\}.$$
Then there exists a unique bijection 
$$\begin{array}{ccc}
 \Irr{H_{\mathbb{C}}} & \to & \mathcal{B} \\
 M &  \mapsto  & V_M
\end{array}$$
such that the following two conditions holds:
\begin{enumerate}
\item For all $V_M \in{\mathcal{B}}$, we have $d_{V_M,M}=1$ and $a({V_M})=a(M)$.
\item If $V\in{\Irr{H_K}}$ and $M\in{\Irr{H_{\mathbb{C}}}}$ are such that $d_{V,M}\neq 0$, then we have $a(M)\leq a(V)$, with equality only for $V=V_M$.
\end{enumerate}
The set $\mathcal{B}$ is called the canonical basic set with respect to the specialization $\theta$.
\end{theorem}
Note that a description of the set $\mathcal{B}$ would lead to a natural parametrization of the set of simple $H_{\mathbb{C}}$-modules. If $H_{\mathbb{C}}$ is semisimple, we know by Tits deformation theorem that the decomposition matrix is just the identity. Hence,  we obtain the following result.
 
\begin{proposition}\label{semisimple}
Assume that $\theta$ is such that $H_{\mathbb{C}}$ is a split semisimple algebra. Then, we have:
$$\mathcal{B}=\Irr{H_K}.$$
\end{proposition}
Now, We want to give an explicit description of $\mathcal{B}$ in the non semisimple case. By \cite[Theorem 7.4.7]{GP}, $H_{\mathbb{C}}$ is semisimple unless $\theta (u)$ is a root of unity. Thus, we can restrict ourselves to the case where $\theta (u)$ is a root of unity. 

The idea is to use results of Ariki which give an interpretation of the decomposition map using the theory of canonical basis of Fock spaces. This theory can not be applied to all Hecke algebras. In fact, this is concerned with the class of Ariki-Koike algebras which contains Hecke algebras of type  $A_{n-1}$  and $B_n$ as special cases.

\section{Ariki-Koike algebras}\label{AK}

a) First, we recall the definition of Ariki-Koike algebras (see \cite{Ma} for a complete  survey of the representation theory of these algebras). Let  $R$ be a commutative ring, let $d\in{\mathbb{N}_{>0}}$,  $n\in{\mathbb{N}}$ and let  $v$, $u_0$, $u_1$,..., $u_{d-1}$ be  $d+1$ parameters in $R$. The Ariki-Koike algebra $\mathcal{H}_{R,n}:=\mathcal{H}_{R,n}(v;u_0,...,u_{d-1})$ (or cyclotomic Hecke algebra of type $G(d,1,n)$) over $R$  is the unital associative $R$-algebra with presentation by:
\begin{itemize}
\item generators: $T_0$, $T_1$,..., $T_{n-1}$,
\item braid relations symbolised by the following diagram:
\\
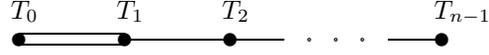
\begin{center}
\begin{picture}(240,20)
\put( 50,10){\circle*{5}}
\put( 47,18){$T_0$}
\put( 50,8){\line(1,0){40}}
\put( 50,12){\line(1,0){40}}
\put( 90,10){\circle*{5}}
\put( 87,18){$T_1$}
\put( 90,10){\line(1,0){40}}
\put(130,10){\circle*{5}}
\put(127,18){$T_2$}
\put(130,10){\line(1,0){20}}
\put(160,10){\circle{1}}
\put(170,10){\circle{1}}
\put(180,10){\circle{1}}
\put(190,10){\line(1,0){20}}
\put(210,10){\circle*{5}}
\put(207,18){$T_{n-1}$}
\end{picture}
\end{center}
and the following ones:
 \begin{align*}
                &(T_0-u_0)(T_0-u_1)...(T_0-u_{d-1})  =  0,\\
                             &(T_i-v)(T_i+1)  =  0\ (i\geq{1}). \end{align*} 
\end{itemize}

These relations are obtained by deforming the relations of the wreath product $(\mathbb{Z}/d\mathbb{Z})\wr{\mathfrak{S}_n}$. We have the following special cases: \begin{itemize} \item if
$d=1$, $\mathcal{H}_{R,n}$ is the Hecke algebra of type $A_{n-1}$ over $R$,  
\item if $d=2$, $\mathcal{H}_{R,n}$ is the Hecke algebra of type $B_{n}$ over $R$.
 \end{itemize}
 It is known that the simple modules of $(\mathbb{Z}/d\mathbb{Z})\wr{\mathfrak{S}_n}$   are indexed by the $d$-tuples of partitions. The same is true for the semisimple Ariki-Koike algebras defined over a field. We say that  $\underline{\lambda}$  is a $d$-partition of rank $n$ if:
\begin{itemize}
\item  $\underline{\lambda}=(\lambda^{(0)},..,\lambda^{(d-1)})$ where, for $i=0,...,d-1$, $\lambda^{(i)}=(\lambda^{(i)}_1,...,\lambda^{(i)}_{r_i})$ is a  partition of rank $|\lambda^{(i)}|$  such that $\lambda^{(i)}_1\geq{...}\geq{ \lambda^{(i)}_{r_i}}>0$,
\item $\displaystyle{\sum_{k=0}^{d-1}{|\lambda^{(k)}|}}=n$.
\end{itemize}
We denote by $\Pi_{n}^d$ the set of  $d$-partitions of rank $n$.

For each $d$-partition $\underline{\lambda}$ of rank $n$, we can associate a $\mathcal{H}_{R,n}$-module $S_R^{\underline{\lambda}}$ which is free over $R$. This is called  a Specht module\footnote{ Here, we use the definition of the classical Specht modules. The passage from classical Specht modules to their duals is provided by the map $(\lambda^{(0)},\lambda^{(1)},...,\lambda^{(d-1)})\mapsto{(\lambda^{(d-1)'},\lambda^{(d-2)'},...,\lambda^{(0)'})}$ where, for $i=0,...,d-1$,  $\lambda^{(i)'}$ denotes the conjugate partition.}. These modules generalize the Specht modules previously defined for Hecke algebras of type $A_{n-1}$ and $B_n$.

Assume that $R$ is a field.  Then, for each $d$-partition of rank $n$, there is a natural bilinear form which is defined over each $S_R^{\underline{\lambda}}$. We denote by $\textrm{rad}$ the radical associated to this bilinear form. 
Then, the non zero  $D_R^{\underline{\lambda}}:=S_R^{\underline{\lambda}}/\textrm{rad}{(S_R^{\underline{\lambda}})}$ form a complete set of non-isomorphic simple $\mathcal{H}_{R,n}$-modules (see for example \cite[chapter 13]{Ab}). In particular, if $\mathcal{H}_{R,n}$ is semisimple, we have $\textrm{rad}(S_R^{\underline{\lambda}})=0$ for all $\underline{\lambda}\in\Pi_{n}^d$ and  the set of simple modules are given by the $S_R^{\underline{\lambda}}$. We have the following  criterion of semisimplicity:

\begin{theorem}[Ariki \cite{A1}]  $\mathcal{H}_{R,n}$ is split semisimple if and only if we have:
\begin{itemize}
\item for all $i\neq{j}$ and for all   $l\in{\mathbb{Z}}$ such that $|l|<{n}$, we have:
$$v^l u_i\neq{u_j},$$
\item $\displaystyle{\prod_{i=1}^{n}{(1+v+...+v^{i-1})}\neq{0}}.$
\end{itemize}
\end{theorem}

Assume that $R$ is a field of characteristic $0$,  using results of Dipper-Mathas \cite{DM}, Ariki-Mathas \cite{AM} and Mathas \cite{Ms},  the case where $\mathcal{H}_{R,n}$ is  not semisimple  can be  reduced to the case where $R=\mathbb{C}$ and where all the $u_i$ are powers of $v$. In this paper, we will mostly concentrate upon the case where $v$ is a primitive root of unity of order $e$.

b) Let $\eta_e:=\textrm{exp} (\frac{2i\pi}{e})$ and let $v_0$, $v_1$,..., $v_{d-1}$ be integers such that $0\leq{v_0}\leq ...\leq v_{d-1}<e$. We consider the Ariki-Koike algebra $\mathcal{H}_{\mathbb{C},n}$  over $\mathbb{C}$ with the following choice of parameters:
\begin{eqnarray*}
&&u_j=\eta_e^{v_{j}}\ \textrm{for}\ j=0,...,d-1,\\
&&v=\eta_{e}. 
\end{eqnarray*}
This algebra is not semisimple in general. Hence the simple $\mathcal{H}_{\mathbb{C},n}$-modules  are given by the non zero $D^{\underline{\lambda}}_{\mathbb{C}}$. We denote:
$$\Gamma^n_0:=\{\underline{\lambda}\in{\Pi_n^d}\ |\ D^{\underline{\lambda}}_{\mathbb{C}}\neq 0 \}.$$
Now,  we wish  to describe the notion of decomposition map in the  context of Ariki-Koike algebras.

Let $y$ be an indeterminate and let $\mathfrak{v}:=y^d$. Let $A=\mathbb{C}[y,y^{-1}]$. We assume that we have $d$ invertible elements $\mathfrak{u}_0,...,\mathfrak{u}_{d-1}$ in $A$ such that we have for all $i\neq j$ and $l\in{\mathbb{Z}}$ with $|l|\leq n$:
$$\mathfrak{v}^l \mathfrak{u}_i \neq \mathfrak{u}_j$$
We consider the Ariki-Koike algebra $\mathcal{H}_{A,n}$ with the following choice of parameters: 
\begin{eqnarray*}
&&u_j=\mathfrak{u}_j\ \textrm{for}\ j=0,...,d-1,\\
&&v=\mathfrak{v}. 
\end{eqnarray*}

 Let $K$ be the field of fractions of $A$ and let $\mathcal{H}_{K,n}:=K\otimes_A\mathcal{H}_{A,n}$. By the above criterion of semisimplicity,  $\mathcal{H}_{K,n}$ is semisimple and its simple modules are given by the Specht modules $S^{\underline{\lambda}}_K$ defined over $K$. Now, let  $\theta :A \to \mathbb{C}$ such that $\theta(\mathfrak{u})=\eta_e:=\textrm{exp}(\frac{2i\pi}{e})$. Assume that we have  $\theta(\mathfrak{u}_j)=\eta_e^{v_j}$ for $j=0,...,d-1$. Then, we have $\mathcal{H}_{\mathbb{C},n}=\mathbb{C}\otimes_A \mathcal{H}_{A,n}$ and the decomposition map is  defined as follows:
$$\begin{array}{cccc}
d_{\theta} : & R_0 (\mathcal{H}_{K,n}) & \to &  R_0 (\mathcal{H}_{\mathbb{C},n})  \\
 & [S_K^{\underline{\lambda}}] & \mapsto &   \displaystyle{[S_{\mathbb{C}}^{\underline{\lambda}}]=\sum_{\underline{\mu}\in{\Gamma^n_0}} d_{ S_K^{\underline{\lambda}},  D^{\underline{\mu}}_{\mathbb{C}}} [D^{\underline{\mu}}_{\mathbb{C}}] }
\end{array}$$
 We will now explain the connections between this decomposition map and the theory of Fock spaces.

\section{Quantum groups and Fock spaces}
 The aim of this part is to introduce Ariki's theorem which provides a way to compute the decomposition maps for Ariki-Koike algebras. For details, we refer to \cite{Ab}. We keep the notations of the previous section.

\subsection{Quantum group of type $A^{(1)}_{e-1}$} Let  $\mathfrak{h}$ be a free $\mathbb{Z}$-module with basis $\{h_i,\mathfrak{d}\ |\ 0\leq{i}<e\}$ and let $\{\Lambda_i,\delta\ |\ 0\leq{i}<e\}$ be the dual basis  with respect to the pairing: 
$$\langle\ ,\ \rangle :\mathfrak{h}^*\times{ \mathfrak{h}}\to{\mathbb{Z}}$$
such that $\langle\Lambda_i,h_j\rangle=\delta_{ij}$, $\langle\delta,\mathfrak{d}\rangle=1$ and  $\langle\Lambda_i,\mathfrak{d}\rangle=\langle\delta,h_j\rangle=0$ for $0\leq{i,j}<e$. For $0\leq{i}<e$, we define the simple roots of $\mathfrak{h}^*$ by:
$$\alpha_i=\left\{ \begin{array}{ll} 2\Lambda_0-\Lambda_{e-1}-\Lambda_1+\delta & \textrm{if}\  i=0, \\
                                   2\Lambda_i-\Lambda_{i-1}-\Lambda_{i+1} & \textrm{if}\  i>0, \\
    \end{array} \right.$$
where $\Lambda_e:=\Lambda_0$. The  $\Lambda_i$ are called  the fundamental weights.

Let $q$ be an indeterminate and let $\mathcal{U}_q$ be the quantum group of type $A^{(1)}_{e-1}$.
This is a unital associative algebra over  $\mathbb{C}(q)$ which is generated by elements $\{e_i,f_i\ |\ i\in{\{0,...,e-1\}}\}$ and $\{k_h\ |\ h\in{\mathfrak{h}}\}$ (see for example \cite{Ab} for the relations).

For $j\in{\mathbb{N}}$ and $l\in{\mathbb{N}}$, we define:
\begin{itemize}
\item $\displaystyle{[j]_q:=\frac{q^j-q^{-j}}{q-q^{-1}}}$,
\item $\displaystyle{ [j]_q^{!}:=[1]_q[2]_q...[j]_q}$,
\item  $\displaystyle{   \left[ \begin{array}{c} l      \\  j \end{array} \right]_q=\frac{[l]_q^!}{ [j]_q^![l-j]_q^!}}$.
\end{itemize}

 Let $\mathcal{A}=\mathbb{Z}[q,q^{-1}]$. We consider the Kostant-Lusztig   $\mathcal{A}$-form of $\mathcal{U}_q$  which is denoted by  $\mathcal{U}_{\mathcal{A}}$: this is a  $\mathcal{A}$-subalgebra of $\mathcal{U}_q$ generated by the divided powers 
$e_i^{(l)}:=\displaystyle{\frac{e^{l}_i}{[l]^{!}_q}} $, $f_j^{(l)}:=\displaystyle{\frac{f^{l}_j}{[l]^{!}_q} }$  for  $0\leq{i,j}<e$ and $l\in{\mathbb{N}}$ and by  $k_{h_{i}}$, $k_{\mathfrak{d}}$,  $k^{-1}_{h_i}$, $k^{-1}_{\mathfrak{d}}$ for $0\leq{i}<e$. Now, if $S$ is a ring and $u$   an invertible element in $S$, we can form the specialized algebra $\mathcal{U}_{S,u}:=S\otimes_{\mathcal{A}}\mathcal{U}_{\mathcal{A}}$  by specializing the indeterminate $q$ to $u\in{S}$.

For any $n\geq 0$, let $\mathcal{F}_n$ be the  $\mathbb{C}(q)$ vector space with basis consisting of all the $d$-partitions of rank $n$. The Fock space is the direct sum:
$$\mathcal{F}:=\oplus_{n\in{\mathbb{N}}} \mathcal{F}_n.$$
We will now see that this space can be endowed with two different structures of  $\mathcal{U}_q$-module. To describe these actions, we need some combinatorial definitions.

Let  $\underline\lambda={(\lambda^{(0)},...,\lambda^{(d-1)})}$ be a $d$-partition of rank  $n$. The diagram of  $\underline{\lambda}$ is the following set:
$$[\underline{\lambda}]=\left\{ (a,b,c)\ |\ 0\leq{c}\leq{d-1},\ 1\leq{b}\leq{\lambda_a^{(c)}}\right\}.$$

The elements of this diagram are called   the nodes of  $\underline{\lambda}$. 
Let  $\gamma=(a,b,c)$ be a node of  $\underline{\lambda}$. The residue of  $\gamma$  associated to the set  $\{e;{v_0},...,{v_{d-1}}\}$ is the element of $\mathbb{Z}/e\mathbb{Z}$ defined by:
$$\textrm{res}{(\gamma)}=(b-a+v_{c})(\textrm{mod}\ e).$$

If $\gamma $ is a node with residue $i$, we say that  $\gamma$ is an  $i$-node. 
Let  $\underline{\lambda}$ and  $\underline{\mu}$ be two $d$-partitions of rank  $n$ and $n+1$ such that  $[\underline{\lambda}]\subset{[\underline{\mu}]}$. There exists a node $\gamma$ such that  $[\underline{\mu}]=[\underline{\lambda}]\cup{\{\gamma\}}$. Then, we denote $[\underline{\mu}]/[\underline{\lambda}]=\gamma$. If $\textrm{res}{(\gamma)}=i$, we say that  $\gamma$ is an addable $i$-node  for   $\underline{\lambda}$ and a removable  $i$-node for  $\underline{\mu}$.

\subsection{Hayachi realization of Fock spaces}

In this part, we consider  the following order on the set of removable and addable nodes of a $d$-partition:  
 we say that $\gamma=(a,b,c)$ is below  $\gamma'=(a',b',c')$ if  $c<c'$ or if  $c=c'$ and $a<a'$.

This order will be called the AM-order and the notion of normal nodes and good nodes below are linked with this order (in the next paragraph, we will give another order on the set of nodes which is distinct from  this one).

 Let $\underline{\lambda}$ and $\underline{\mu}$ be two $d$-partitions of rank $n$ and $n+1$ such that  there exists an $i$-node $\gamma$ such that  $[\underline{\mu}]=[\underline{\lambda}]\cup{\{\gamma\}}$. We define:
\begin{align*}
N_i^{a}{(\underline{\lambda},\underline{\mu})}=&   \sharp\{ \textrm{addable }\ i-\textrm{nodes of } \underline{\lambda}\ \textrm{ above } \gamma\} \\
 & -\sharp\{ \textrm{removable }\ i-\textrm{nodes of } \underline{\mu}\ \textrm{ above } \gamma\},\\
N_i^{b}{(\underline{\lambda},\underline{\mu})}= &   \sharp\{ \textrm{addable } i-\textrm{nodes of } \underline{\lambda}\ \textrm{ below } \gamma\}\\
&   -\sharp\{ \textrm{removable } i-\textrm{nodes of } \underline{\mu}\ \textrm{ below } \gamma\},\\
N_{i}{(\underline{\lambda})} =& \sharp\{ \textrm{addable } i-\textrm{nodes of } \underline{\lambda}\}\\
& -\sharp\{ \textrm{removable } i-\textrm{nodes of } \underline{\lambda}\},\\
N_{\mathfrak{d}}{(\underline{\lambda})} =& \sharp\{ 0-\textrm{nodes of } \underline{\lambda}\}.
\end{align*}

\begin{theorem}[Hayashi \cite{Ha}]  $\mathcal{F}$ becomes a $\mathcal{U}_q$-module with respect to the following action:
$$e_{i}\underline{\lambda}=\sum_{\textrm{res}([\underline{\lambda}]/[\underline{\mu}])=i}{q^{-N_i^{a}{(\underline{\mu},\underline{\lambda})}}\underline{\mu}},\qquad{f_{i}\underline{\lambda}=\sum_{\textrm{res}([\underline{\mu}]/[\underline{\lambda}])=i}{q^{N_i^{b}{(\underline{\lambda},\underline{\mu})}}\underline{\mu}}},$$
$$k_{h_i}\underline{\lambda}=q^{N_i{(\underline{\lambda})}}\underline{\lambda},\qquad{k_{\mathfrak{d}}\underline{\lambda}=q^{-N_{\mathfrak{d}}{(\underline{\lambda})}}\underline{\lambda}},$$
where $i=0,...,e-1$.
\end{theorem}

Let $\mathcal{M}$ be the $\mathcal{U}_q$-submodule of  $\mathcal{F}$ generated by the empty $d$-partition. It is isomorphic to an integrable  highest weight module. In \cite{K} and \cite{L3},  Kashiwara and Lusztig have independantly shown the existence of a remarkable basis for this class of modules: the  canonical basis. We will see the links between the canonical basis of  $\mathcal{M}$ and the decomposition map for $\mathcal{H}_{\mathbb{C},n}$. First, it is known that the elements of this basis are labeled by the vertices of a certain graph called the crystal graph. 

Based on Misra and Miwa's result,  Ariki and Mathas observed that the vertices of this graph are given by the set of Kleshchev $d$-partitions which we will now define.

Let   $\underline{\lambda}$ be a   $d$-partition and  let $\gamma$ be  an $i$-node, we say that  $\gamma$ is a normal $i$-node of  $\underline{\lambda}$ if, whenever $\eta$ is an $i$-node of $\underline{\lambda}$   below  $\gamma$, there are more removable $i$-nodes between $\eta$ and $\gamma$ than addable $i$-nodes between  $\eta$ and $\gamma$. 
If $\gamma$ is the highest normal $i$-node of    $\underline{\lambda}$, we say that  $\gamma$ is a good $i$-node.

We can now define the notion of Kleshchev $d$-partitions  associated to the set  $\{e;{v_0},...,{v_{d-1}}\}$:
\begin{definition}
The Kleshchev $d$-partitions are defined   recursively  as follows.
\begin{itemize}
   \item The empty partition $\underline{\emptyset}:=(\emptyset,\emptyset,...,\emptyset)$ is  Kleshchev.
    \item If $\underline{\lambda}$ is  Kleshchev, there exist $i\in{\{0,...,e-1\}}$ and a good $i$-node $\gamma$ such that if we remove $\gamma$ from  $\underline{\lambda}$, the resulting  $d$-partition is Kleshchev.
\end{itemize}
\end{definition}
We denote  by $\Lambda^{0,n}_{\{e;v_0,...,v_{l-1}\}}$ the set of  Kleshchev $d$-partitions associated to  the set  $\{e;{v_0},...,{v_{d-1}}\}$. If there is no ambiguity concerning  $\{e;{v_0},...,{v_{d-1}}\}$, we denote it by $\Lambda^{0,n}$.
Now, the crystal graph of $\mathcal{M}$ is given by:
\begin{itemize}
\item vertices: the  Kleshchev $d$-partitions,
\item edges: $\displaystyle{{\underline{\lambda}\overset{i}{\rightarrow}{{\underline{\mu}}}}}$ if and only if $[\underline{\mu}]/[\underline{\lambda}]$ is a good  $i$-node.
\end{itemize}
Thus, the canonical basis $\mathfrak{B}$ of $\mathcal{M}$ is labeled by  the Kleshchev $d$-partitions:

$$\mathfrak{B}=\{ G(\underline{\lambda})\ |\ \underline{\lambda}\in{\Lambda^0_{\{e;v_0,...,v_{d-1}\}}},\ n\in{\mathbb{N}}\}.$$
 
This set is a basis of  the $\mathcal{U}_{\mathcal{A}}$-module $\mathcal{M}_{\mathcal{A}}$ generated by the empty $d$-partition and for any specialization  of $q$ into an invertible element ${u}$  of  a field $R$, we obtain a basis of the specialized module  $\mathcal{M}_{R,u}$   by specializing the set  $\mathfrak{B}$. 

By the characterization of the canonical basis,  for each $\underline{\lambda}\in{\Lambda^{0,n}}$, there exist polynomials $d_{\underline{\mu},\underline{\lambda}}(q)\in{\mathbb{Z}[q]}$ and   a unique element  $G(\underline{\lambda})$ of the canonical basis such that:
$$G(\underline{\lambda})=\sum_{\underline{\mu}\in{\Pi_n^d}}{d_{\underline{\mu},\underline{\lambda}}(q)\underline{\mu}}\qquad{\textrm{and}}\qquad{G(\underline{\lambda})=\underline{\lambda}\ (\textrm{mod}\ q)}.$$

 Now, we have the following theorem of Ariki which shows that the problem of computing the decomposition numbers of  $\mathcal{H}_{R,n}$ can be translated to that of computing the canonical basis of $\mathcal{M}$.
 This theorem was first conjectured by Lascoux, Leclerc and Thibon (\cite{LLT}) in the case of  Hecke algebras of type $A_{n-1}$.

\begin{theorem}[Ariki \cite{A2}]\label{Ar}  There exist a bijection $j_0 : \Lambda^{0,n} \to \Gamma^n_{0}$ such that 
for all  $\underline{\mu}\in{\Pi_n^d}$ and $\underline{\lambda}\in{\Lambda^{0,n}}$, we have:
  $$d_{\underline{\mu},\underline{\lambda}}(1)=d_{S_K^{\underline{\mu}},D_{\mathbb{C}}^{j_0(\underline{\lambda})}},$$
where we recall that:
$$\Gamma^n_0:=\{\underline{\lambda}\in{\Pi_n^d}\ |\ D^{\underline{\lambda}}_{\mathbb{C}}\neq 0 \}.$$

\end{theorem}
Moreover, we have:
\begin{theorem}[Ariki \cite{A3}, Ariki-Mathas \cite{AM}]\label{Ar2}
We have $\Gamma^n_0=\Lambda^{0,n}$ and $j_0=Id$.
\end{theorem}

Hence the above theorem gives a first classification of the simple modules by the set of Kleshchev $d$-partitions. As noted in the introduction,   the problem of this parametrization of the simple $\mathcal{H}_{\mathbb{C},n}$-modules  is that we only know  a recursive description of the Kleshchev $d$-partitions.  We  now deal with another parametrization of this set found by Foda et  al. which uses almost the same objects as Ariki and Mathas.

 \subsection{JMMO realization of Fock space} This action has been defined in \cite{JMMO} and has been used and studied in \cite{FLOTW}. We need to define another order on the set of nodes of a $d$-partitions. 

 Here, we say that  $\gamma=(a,b,c)$ is above  $\gamma'=(a',b',c')$ if:
$$b-a+v_c<b'-a'+v_{c'}\ \textrm{or } \textrm{if}\ b-a+v_c=b'-a'+v_{c'}\textrm{ and }c>c'.$$
This order will be called the FLOTW order and it  allows us to define functions $\overline{N}_i^{a}{(\underline{\lambda},\underline{\mu})}$ and  $\overline{N}_i^{b}{(\underline{\lambda},\underline{\mu})}$ given by the same way as ${N}_i^{a}{(\underline{\lambda},\underline{\mu})}$ et     ${N}_i^{b}{(\underline{\lambda},\underline{\mu})}$ for the AM order.

Now, we have the following result:

\begin{theorem}[Jimbo, Misra, Miwa, Okado \cite{JMMO}] $\mathcal{F}$ is a  $\mathcal{U}_q$-module with respect to the action:
$$e_{i}\underline{\lambda}=\sum_{\operatorname{res}([\underline{\lambda}]/[\underline{\mu}])=i}{q^{-\overline{N}_i^{a}{(\underline{\mu},\underline{\lambda})}}\underline{\mu}},\qquad{f_{i}\underline{\lambda}=\sum_{\operatorname{res}([\underline{\mu}]/[\underline{\lambda}])=i}{q^{\overline{N}_i^{b}{(\underline{\lambda},\underline{\mu})}}\underline{\mu}}},$$
$$k_{h_i}\lambda=v^{{N}_i{(\underline{\lambda})}}\underline{\lambda},\qquad{k_{\mathfrak{d}}\underline{\lambda}=q^{-N_{\mathfrak{d}}{(\underline{\lambda})}}\underline{\lambda}},$$
 where $0\leq{i}\leq{n-1}$. This action will be called the JMMO action.
\end{theorem}

We denote by $\overline{\mathcal{M}}$ the $\mathcal{U}_q$-module generated by the empty $d$-partition with respect to the above action. This is a highest weight module which is isomorphic to  $\mathcal{M}$. However, the elements of the canonical basis are differents in general. Here,  the  $d$-partitions of the crystal graph are obtained  recursively by adding good nodes to  $d$-partitions of the crystal graph with respect to the FLOTW order. 

Foda et al. showed that the analogue of the notion of Kleshchev $d$-partitions for this action is as follows: 

\begin{definition}[Foda, Leclerc, Okado, Thibon, Welsh \cite{FLOTW}] 
 We say that  $\underline\lambda={(\lambda^{(0)} ,...,\lambda^{(d-1)})}$
is a FLOTW  $d$-partition associated to  the set  ${\{e;{v_0},...,v_{d-1}\}}$ if and only if:
\begin{enumerate}
\item for all $0\leq{j}\leq{d-2}$ and $i=1,2,...$, we have:
\begin{align*}
&\lambda_i^{(j)}\geq{\lambda^{(j+1)}_{i+v_{j+1}-v_j}},\\
&\lambda^{(d-1)}_i\geq{\lambda^{(0)}_{i+e+v_0-v_{d-1}}};
\end{align*}
\item  for all  $k>0$, among the residues appearing at the right ends of the length $k$ rows of   $\underline\lambda$, at least one element of  $\{0,1,...,e-1\}$ does not occur.
\end{enumerate}
We denote by $\Lambda^{1,n}_{\{e;{v_0},...,{v_{d-1}}\}}$ the set of FLOTW $d$-partitions of rank $n$ associated to the set ${\{e;{v_0},...,{v_{d-1}}\}}$. If there is no ambiguity concerning  $\{e;{v_0},...,{v_{d-1}}\}$, we denote it by $\Lambda^{1,n}$.
\end{definition}

Hence, the crystal graph of   $\overline{\mathcal{M}}$ is given by:
\begin{itemize}
\item   vertices: the FLOTW $d$-partitions,
\item  edges: $\displaystyle{{\underline{\lambda}\overset{i}{\rightarrow}{{\underline{\mu}}}}}$ if and only if  $[\underline{\mu}]/[\underline{\lambda}]$ is  good  $i$-node with respect  to the FLOTW order.
\end{itemize}

So, the canonical basis elements of $\overline{\mathcal{M}}$ are labeled by the  FLOTW $d$-partitions:
$$\overline{\mathfrak{B}}=\{ \overline{G}(\underline{\lambda})\ |\ \underline{\lambda}\in{\Lambda^{1,n}_{\{e;v_0,...,v_{d-1}\}}},\ n\in{\mathbb{N}}\}.$$
 If we specialize these elements to $q=1$, we obtain the same elements as  in Theorem \ref{Ar} (note that the action of the quantum group on the Fock space specialized at $q=1$ leads to the same module stucture for the Hayashi action and for the JMMO action).

By the characterization of the canonical basis,  for each $\underline{\lambda}\in{\Lambda^{1,n}}$, there exist polynomials $\overline{d}_{\underline{\mu},\underline{\lambda}}(q)\in{\mathbb{Z}[q]}$ and   a unique element  $\overline{G}(\underline{\lambda})$ of the canonical basis such that:
$$\overline{G}(\underline{\lambda})=\sum_{\underline{\mu}\in{\Pi_n^d}}{\overline{d}_{\underline{\mu},\underline{\lambda}}(q)\underline{\mu}}\qquad{\textrm{and}}\qquad{\overline{G}(\underline{\lambda})=\underline{\lambda}\ (\textrm{mod}\ q)}.$$

By Ariki's theorem, we have:
\begin{theorem}[Ariki \cite{A2}]\label{AR2}
There exist a bijection $j_1 : \Lambda^{1,n} \to \Gamma_0^n=\Lambda^{0,n}$ such that 
for all  $\underline{\mu}\in{\Pi_n^d}$ and $\underline{\lambda}\in{\Lambda^{0,n}}$, we have:
  $$\overline{d}_{\underline{\mu},\underline{\lambda}}(1)=d_{S_K^{\underline{\mu}},D_{\mathbb{C}}^{j_1(\underline{\lambda})}}.$$
\end{theorem}
Hence, we can alternatively use the JMMO action instead of  the Hayashi action to compute the decomposition matrix for Ariki-Koike algebras. A natural question is now to ask if there is an  interpretation of the FLOTW multipartitions in  the representation theory of Ariki-Koike algebras. An answer will be given by extending the results of Geck and Rouquier to the case of Ariki-Koike algebras.

\section{Canonical basic sets for Ariki-Koike algebras}

Let  $e$ be a positive integer, $\eta_e:=\textrm{exp} (\frac{2i\pi}{e})$ and let $v_0$, $v_1$,..., $v_{d-1}$ be integers such that $0\leq{v_0}\leq ...\leq v_{d-1}<e$. We consider the Ariki-Koike algebra $\mathcal{H}_{\mathbb{C},n}$  over $\mathbb{C}$ with the following choice of parameters:
\begin{eqnarray*}
&&u_j=\eta_e^{v_{j}}\ \textrm{for}\ j=0,...,d-1,\\
&&v=\eta_{e}, 
\end{eqnarray*}

 In this part,  we show that there exists a  ``canonical basic set'' of Specht modules which is in bijection with the set of simple   $\mathcal{H}_{\mathbb{C},n}$-modules. To do this, we consider the  Ariki-Koike algebra   $\mathcal{H}_{\mathbb{C},n}$   and we study the Kashiwara-Lusztig  canonical basis of the associated highest weight $\mathcal{U}(\widehat{sl}_e)$-module. First, we have to define a semisimple Ariki-Koike algebra which can be specialized to   $\mathcal{H}_{\mathbb{C},n}$ as in the end of \S \ref{AK} b).  Let $y$ be an indeterminate and let $\mathfrak{v}=y^d$. Let $A=\mathbb{C}[y,y^{-1}]$.

We consider the Ariki-Koike algebra $\mathcal{H}_{A,n}$ with the following choice of parameters: 
\begin{eqnarray*}
&&u_j=\eta_d^j y^{dv_j-je+de}\ \textrm{for}\ j=0,...,d-1,\\
&&v=\mathfrak{v}, 
\end{eqnarray*}
where $\eta_d:=\textrm{exp}(\frac{2i\pi}{d})$. Let $K:=\mathbb{C}(y)$ and  let $\mathcal{H}_{K,n}:=K\otimes_A \mathcal{H}_{A,n}$. It is easy to see that this algebra is semisimple and that, under the specialization $y\in{A}\mapsto \textrm{exp}(\frac{2i\pi}{de})\in{\mathbb{C}}$, we obtain the algebra $\mathcal{H}_{\mathbb{C},n}$. Hence, we have a decomposition map as in \S \ref{AK} b):
$$\begin{array}{cccc}
d_{\theta} : & R_0 (\mathcal{H}_{K,n}) & \to &  R_0 (\mathcal{H}_{\mathbb{C},n})  \\
 & [S_K^{\underline{\lambda}}] & \mapsto &   \displaystyle{[S_{\mathbb{C}}^{\underline{\lambda}}]=\sum_{\underline{\mu}\in{\Gamma^n_0}} d_{ S_K^{\underline{\lambda}},  D^{\underline{\mu}}_{\mathbb{C}}} [D^{\underline{\mu}}_{\mathbb{C}}] }
\end{array}$$

Now, we have to define an $a$-value on the simple modules (which are the Specht modules defined over $K$). The main problem here is that we don't have Kazhdan-Lusztig type bases for  Ariki-Koike algebras in general but we do have Schur elements: they have been computed by Geck, Iancu and Malle in \cite{GIM}. This leads to  the following definition of $a$-values (see \cite[\S 3.3]{J2}):
\begin{definition}\label{afonction} Let $\underline{\lambda}:=(\lambda^{(0)},\lambda^{(1)},...,\lambda^{(d)})\in{\Lambda^{1,n}}$ where for $j=0,...,d-1$ we have $\lambda^{(j)}:=(\lambda^{(j)}_1,...,\lambda^{(j)}_{h^{(j)}})$. We assume that the rank of  $\underline{\lambda}$ is $n$.   For $j=0,...,d-1$ and $p=1,...,n$, we define the following rational numbers:
\begin{align*}
 m^{(j)}&:=v_j-\frac{ie}{d}+e,\\
B'^{(j)}_p&:=\lambda^{(j)}_p-p+n+m^{(j)},\end{align*}
where we use the convention that $\lambda^{(j)}_p:=0$ if $p>h^{(j)}$. For $j=0,...,d-1$, let $B'^{(j)}=(B'^{(j)}_1,...,B'^{(j)}_n)$.   Then, we define:
$$a_1(\underline{\lambda}):=\sum_{{0\leq{i}\leq{j}<d}\atop{{(a,b)\in{B'^{(i)}\times{B'^{(j)}}}}\atop{a>b\ \textrm{if}\ i=j}}}{\min{\{a,b\}}}   -\sum_{{0\leq{i,j}<d}\atop{{a\in{B'^{(i)}}}\atop{1\leq{k}\leq{a}}}}{\min{\{k,m^{(j)}\}}}.$$
\end{definition}
 The $a$-value associated to  $S^{\underline{\lambda}}_K$ is the rational number $a(\underline{\lambda}):=a_1(\underline{\lambda})+f(n)$ where $f(n)$ is a rational number which only depends on the parameters $\{e;v_0,...,v_{d-1}\}$ and on $n$ (the expression of $f$ is given in \cite{J2}).

Next, we associate to  each  $\underline{\lambda}\in{\Lambda^{1,n}}$ a sequence of residues which will have ``nice'' properties with respect to the $a$-value:
\begin{proposition}[\cite{J2}]\label{asuite} Let $\underline{\lambda}\in{\Lambda^{1,n}}$ and let:
$$l_{\textrm{max}}:=\operatorname{max}\{\lambda^{(0)}_1,...,\lambda^{(l-1)}_1\}.$$
Then,  there exists a removable node $\xi_1$ with residue  $k$ on a part  $\lambda_{j_1}^{(i_1)}$ with  length $l_{\textrm{max}}$, such that there doesn't exist a $k-1$-node at the right end of a part with length $l_{\textrm{max}}$ (the existence of such a node is proved in \cite[Lemma 4.2]{J2}). 

Let $\gamma_1$, $\gamma_2$,..., $\gamma_r$ be the  $k-1$-nodes at the right ends of parts $\displaystyle {\lambda_{p_1}^{(l_1)}\geq{\lambda_{p_2}^{(l_2)}\geq{...}}}$$\geq{\lambda_{p_r}^{(l_r)}}$. 
 Let   $\xi_1$, $\xi_2$,..., $\xi_{s}$ be the  removable $k$-nodes  of   $\underline{\lambda}$ on  parts  $\displaystyle{\lambda_{j_1}^{(i_1)}\geq{\lambda_{j_2}^{(i_2)}\geq{...}\geq{\lambda_{j_s}^{(i_s)}}}}$   such that:
$$\lambda_{j_s}^{(i_s)}>\lambda_{p_1}^{(l_1)}.$$

We remove the nodes  $\xi_1$, $\xi_2$,..., $\xi_s$ from $\underline{\lambda}$. Let  $\underline{\lambda}'$ be the  resulting $d$-partition. Then,   $\underline{\lambda}'\in{\Lambda}^{1,n-s}$ and we define recursively the  $a$-sequence  of residues of  $\underline{\lambda}$ by:
$$a\textrm{-sequence}(\underline{\lambda})=a\textrm{-sequence}(\underline{\lambda}'),\underbrace{k,...,k}_{s}.$$
\end{proposition}
$\ $\\
\textbf{Example}:\\
Let $e=4$, $d=3$, $v_0=0$, $v_1=2$ and $v_2=3$. We consider the  $3$-partition $\underline{\lambda}=(1,3.1,2.1.1)$ with the following diagram:
$$  \left( \ \begin{tabular}{|c|}
                         \hline
                           0    \\
                          \hline
                             \end{tabular}\  ,\
                       \begin{tabular}{|c|c|c|}
                         \hline
                           2  & 3  & 0 \\
                          \hline
                          1     \\
                          \cline{1-1}

                             \end{tabular}\ ,\
                        \begin{tabular}{|c|c|}
                         \hline
                           3  & 0 \\
                          \hline
                          2     \\
                          \cline{1-1}
                          1 \\
                           \cline{1-1}
                             \end{tabular}\
                                                 \right)$$
\\
$\underline{\lambda}$ is a FLOTW $3$-partition.

We want to determine the $a$-sequence of $\underline{\lambda}$: we have to find  $k\in{\{0,1,2,3\}}$, $s\in{\mathbb{N}}$ and a $2$-partition $\underline{\lambda}'$ such that:
$$a\textrm{-sequence}(\underline{\lambda})=a\textrm{-sequence}(\underline{\lambda}'),\underbrace{k,...,k}_{s}.$$
The part with maximal length is  the  part with length  $3$ and the residue of the associated removable node is $0$. We remark that there are two others removable  $0$-nodes on parts with length $1$ and $2$. Since there is no  node with residue $0-1\equiv 3\ (\textrm{mod}\ e)$ at the right ends of the parts of  $\underline{\lambda}$, we must remove these three $0$-nodes. Thus, we have to take   $k=0$, $s=3$ and   $\underline{\lambda}'=(\emptyset,2.1,1.1.1)$, hence:
$$a\textrm{-sequence}(\underline{\lambda})=a\textrm{-sequence}(\emptyset,2.1,1.1.1),0,0,0.$$
Observe  that the  $3$-partition $(\emptyset,2.1,1.1.1)$ is a FLOTW  $3$-partition.

Now, the residue of the removable node on the part with maximal length is $3$. Thus, we obtain:
$$a\textrm{-sequence}(\underline{\lambda})=a\textrm{-sequence}(\emptyset,1.1,1.1.1),3,0,0,0.$$
Repeating the same procedure, we finally  obtain:
$$a\textrm{-sequence}(\underline{\lambda})=3,2,2,1,1,3,0,0,0.$$
\\
\begin{proposition}[\cite{J2}]\label{decompo} Let $n\in{\mathbb{N}}$, let $\underline{\lambda}\in{\Lambda}^{1,n}$  and let $a\textrm{-sequence}(\underline{\lambda})=\underbrace{i_1,...,i_1}_{a_1},\underbrace{i_{2},...,i_{2}}_{a_{2}},...,\underbrace{i_s,...,i_s}_{a_s}$ be its $a$-sequence of residues where we assume that for all $j=1,...,s-1$, we have $i_{j}\neq{i_{j+1}}$. Then, we have:

$$A(\underline{\lambda}):=f^{(a_s)}_{i_s}f^{(a_{s-1})}_{i_{s-1}} ...f^{(a_1)}_{i_1}\underline{\emptyset}=\underline{\lambda}+ \sum_{a(\underline{\mu})>a(\underline{\lambda})}{c_{\underline{\lambda},\underline{\mu}}(q)\underline{\mu}},$$
 where $c_{\underline{\lambda},\underline{\mu}}(q)\in{\mathbb{Z}[q,q^{-1}]}$.
\end{proposition}
It is obvious that the set $\{A(\underline{\lambda})\ |\ \underline{\lambda}\in{\Lambda}^{1,n},\  n\in{\mathbb{N}}\}$ is a basis of ${\mathcal{M}}_{\mathcal{A}}$. Using the characterization of the canonical basis, we obtain the following theorem:
\begin{proposition}[\cite{J2}]\label{base} Let $n\in{\mathbb{N}}$ and let $\underline{\lambda}\in{\Lambda}^{1,n}$, then we have:
$$\overline{G}({\underline{\lambda}}):=\underline{\lambda}+\sum_{a(\underline{\mu})>a(\underline{\lambda})}\overline{d}_{\underline{\mu},\underline{\lambda}}(q)\underline{\mu}$$
\end{proposition}
Now, assume that $\underline{\mu}\in{\Lambda^{0,n}}$ is a Kleshchev multipartition and let $G({\underline{\mu}})$ be the element of the canonical basis of $\mathcal{M}$ such that $G({\underline{\mu}})=\underline{\mu}\ (\textrm{mod}\ v)$. Let $\overline{G}(c({\underline{\mu}}))$ with  $c(\underline{\mu})\in{\Lambda^{1,n}}$ be the element of the canonical basis of  $\overline{\mathcal{M}}$ such that  $\overline{G}(c({\underline{\mu}}))$ coincides with  $G({\underline{\mu}})$ at $v=1$. This defines a bijection:
$$\begin{array}{cccc}
c : & \Lambda^{0,n} & \to &  \Lambda^{1,n}  
\end{array}$$
This bijection can be described by reading the crystal graphs of $\mathcal{M}$ and $\overline{\mathcal{M}}$ (see \cite{FLOTW}).

We can now define an $a$-value on the set of simple $\mathcal{H}_{\mathbb{C},n}$-modules by setting for  $\underline{\mu}\in{\Lambda^{0,n}=\Gamma^n_0}$:
$$a(D^{\underline{\mu}}_{\mathbb{C}}):=a(S^{c(\underline{\mu})}_{K}).$$

Combining the above proposition with  Theorem \ref{AR2}, we obtain the following result which shows the existence of the canonical basic set for Ariki-Koike algebras and gives an explicit description of this set. 

\begin{theorem}[\cite{J2}]\label{cbs} 
 We define the following subset of $\Irr{\mathcal{H}_{K,n}}$:
$$\mathcal{B}:=\{  S^{\underline{\lambda}}_{K} \in{\Irr{{\mathcal{H}_{K,n}}}}\ |\ \underline{\lambda}\in{\Lambda^{1,n}}\}.$$
Then:
\begin{enumerate}
\item For all $\underline{\mu}\in{\Lambda^{0,n}}$, we have $d_{S^{c(\underline{\mu})}_{K} ,D^{\underline{\mu}}_{\mathbb{C}}    }=1$.
\item If $V\in{\Irr{\mathcal{H}_{K,n}}}$ and $M\in{\Irr{\mathcal{H}_{\mathbb{C},n}}}$ are such that $d_{V,M}\neq 0$, then we have $a(M) \leq a(V)$, with equality only if $M=D^{\underline{\mu}}_{\mathbb{C}}$ with $\underline{\mu}\in{\Lambda^{0,n}}$ and $V=S^{c(\underline{\mu})}_{K}$ 
\item We have $c=j_1^{-1}$ where $j_1$ is as in Theorem \ref{AR2}.
\end{enumerate}
The set $\mathcal{B}$ is called the canonical basic set with respect to the specialization $\theta$ and it is in natural bijection with $\Irr{\mathcal{H}_{\mathbb{C},n}}$.
\end{theorem}
This theorem shows that the decomposition matrix for Ariki-Koike algebras has a lower triangular shape if we order the rows and columns with respect to the $a$-value.

Note also that the results of this section induce a purely combinatorial algorithm for the computation of the canonical basis and of the decomposition matrices of Ariki-Koike algebras (see \cite{J4}) which generalizes the LLT algorithm.

\section{Consequences}

We obtain an explicit description of the canonical basic set for all Hecke algebras of classical type in characteristic $0$.
  
\subsection{Type $A_{n-1}$} Assume that $W$ is a Weyl group of type $A_{n-1}$ and that $\theta (u)$ is a primitive $e^{\textrm{th}}$-root of unity. $H$ is a special case of Ariki-Koike algebras. Hence, we can use Theorem \ref{cbs} to find the canonical basic sets. We note that we can also find this set using results of Dipper and James \label{DJ} as it is expained in \cite[Example 3.5]{Gk}.

\begin{proposition} Assume that $W$ is a Weyl group of type $A_{n-1}$ and that $\theta (u)$ is a primitive $e^{\textrm{th}}$-root of unity. Then, we have:
$$\mathcal{B}=\{S^{\lambda}_K\ |\ \lambda\in{\Lambda^{1,n}_{\{e;0\}}}\}.$$
Note that:
$$\begin{array}{lll}
 \lambda\in{\Lambda^{1,n}_{\{e;0\}}} &\iff & |\lambda|=n\textrm{ and }\lambda=(\lambda_1,...,\lambda_r)\textrm{  is }e-\textrm{regular}.\\
  &\iff& |\lambda|=n\textrm{ and } \textrm{for all }i\in{\mathbb{N}},\textrm{we can't have  } \lambda_i=...=\lambda_{i+e-1}\neq 0.
\end{array}$$
\end{proposition}

\subsection{Type $B_n$} Assume that $W$ is a Weyl group of type $B_n$ and that $\theta (u)$ is a primitive $e^{\textrm{th}}$-root of unity. $H$ is a special case of Ariki-Koike algebras. Hence, we can use Theorem \ref{cbs} to find the canonical basic sets. 

\begin{proposition} Assume that $W$ is a Weyl group of type $B_n$ and that $\theta (u)$ is a primitive $e^{\textrm{th}}$-root of unity. Then, we have:
\begin{itemize}
\item if $e$ is odd:
$$\mathcal{B}=\{S_K^{(\lambda^{(0)},\lambda^{(1)})}\ |\ \lambda^{(0)}\in{\Lambda^{1,n_0}_{\{e;0\}}},\ \lambda^{(1)}\in{\Lambda^{1,n_1}_{\{e;0\}}},\ n_0+n_1=n\}.$$
\item  if $e$ is even:
$$\mathcal{B}=\{S^{(\lambda^{(0)},\lambda^{(1)})}\ |\ (\lambda^{(0)},\lambda^{(1)})\in{\Lambda^{1,n}_{\{e;1,\frac{e}{2}\}}}\}.$$
Recall that $(\lambda^{(0)},\lambda^{(1)})\in{\Lambda^{1,n}_{\{e;1,\frac{e}{2}\}}}$ if and only if $|\lambda^{(0)}|+|\lambda^{(1)}|=n$ and:
\begin{enumerate}
\item for all  $i=1,2,...$, we have:
\begin{align*}
&\lambda_i^{(0)}\geq{\lambda^{(1)}_{i+\frac{e}{2}-1}},\\
&\lambda^{(1)}_i\geq{\lambda^{(0)}_{i+\frac{e}{2}+1}};
\end{align*}
\item  for all  $k>0$, among the residues appearing at the right ends of the length $k$ rows of   $(\lambda^{(0)},\lambda^{(1)})$, at least one element of  $\{0,1,...,e-1\}$ does not occur.
\end{enumerate}
\end{itemize}
\end{proposition}

\subsection{Type $D_n$} Assume that $W$ is a Weyl group of type $D_n$ and that $\theta (u)$ is a primitive $e^{\textrm{th}}$-root of unity. Then, $H$ can be seen as a subalgebra of an Hecke algebra $H_1$ of type $B_n$ as in \S \ref{dd}. 
The specialization $\theta$ induces a decomposition map for $H_1$:
$$d^1_{\theta}:R_0(H_{1,K})\to R_0(H_{1,\mathbb{C}}).$$
Hecke algebras of type $B_n$ with unequal parameters are special cases of Ariki-Koike algebras. Hence, we can also define a canonical basic set for these algebras (the existence has been previously proved in \cite{Gex}). Furthermore, in \cite{Gex}, Geck has shown that the simple $H_K$-modules in  the canonical  basic set for type $D_n$ are those which appear in the restriction of the simple $H_{1,K}$-modules of the canonical basic set for type $D_n$. We obtain the following description of $\mathcal{B}$.

\begin{proposition} Assume that $W$ is a Weyl group of type $D_n$ and that $\theta (u)$ is a primitive $e^{\textrm{th}}$-root of unity. Then:
\begin{itemize}
\item if $e$ is odd, we have:
$$\begin{array}{c}
\mathcal{B}=\left\{V^{[\lambda^{(0)},\lambda^{(1)}]}\ |\ \lambda^{(0)}\neq{\lambda^{(1)}}, \lambda^{(0)}\in{\Lambda^{1,n_0}_{\{e;0\}}},\ \lambda^{(1)}\in{\Lambda^{1,n_1}_{\{e;0\}}},\ n_0+n_1=n      \right\}\\
\bigcup{\left\{V^{[\lambda^{(0)},\pm]}\ |\ \lambda^{(0)}\in{\Lambda^{1,\frac{n}{2}}_{\{e;0\}}}\right\}}.\end{array}$$
\item  if $e$ is even, we have:
$$\begin{array}{c}
\mathcal{B}=\left\{V^{[\lambda^{(0)},\lambda^{(1)}]}\ |\ \lambda^{(0)}\neq{\lambda^{(1)}},\ (\lambda^{(0)},\lambda^{(1)})\in{\Lambda^{1,n}_{\{e;0,\frac{e}{2}\}}}\right\}\\
\bigcup{\left\{V^{[\lambda^{(0)},\pm]}\ |\ (\lambda^{(0)},\lambda^{(0)})\in{\Lambda^{1,n}_{\{e;0,\frac{e}{2}\}}} \right\}}.\end{array}$$
Recall that $(\lambda^{(0)},\lambda^{(1)})\in{\Lambda^1_{\{e;0,\frac{e}{2}\}}}$ if and only if $|\lambda^{(0)}|+|\lambda^{(1)}|=n$ and:
\begin{enumerate}
\item for all  $i=1,2,...$, we have:
\begin{align*}
&\lambda_i^{(0)}\geq{\lambda^{(1)}_{i+\frac{e}{2}}},\\
&\lambda^{(1)}_i\geq{\lambda^{(0)}_{i+\frac{e}{2}}};
\end{align*}
\item  for all  $k>0$, among the residues appearing at the right ends of the length $k$ rows of   $(\lambda^{(0)},\lambda^{(1)})$, at least one element of  $\{0,1,...,e-1\}$ does not occur.
\end{enumerate}
\end{itemize}
\end{proposition}
Note that Hu  \cite{H} has given another parametrization of the simple ${H}_{\mathbb{C}}$-modules by using the Kleshchev bi-partitions. The connection with the above result is explained in \cite[\S 4.3]{J3}

\subsection{Exceptional types} 
An explicit description of the canonical basic set for all exceptional types and for all specializations can be found in \cite[Chapter 3]{J3}.

\subsection{Positive characteristic} 
The existence of the canonical basic set for Hecke algebras of finite Weyl group $W$ can also be proved when the Hecke algebra is defined over a field of characteristic $p$ where $p$ is a ``good'' prime number for $W$ (see \cite{GR}). In fact, it is easy to see that the parametrization of the canonical basic set in characteristic $0$ holds in ``good'' characteristic (see \cite{J5}).

\subsection{Cyclotomic Hecke algebras of type $G(d,p,n)$} Both Ariki-Koike algebras and Hecke algebras of type $D_n$ are particular cases of cyclotomic Hecke algebras of type $G(d,p,n)$ (with $p|d$). This kind of algebras have been defined in \cite{Ac} and   can be seen as subalgebras of Ariki-Koike algebras. By using Clifford theory and results on graded algebras proved by Genet in \cite{Gg}, we can obtained a parametrization of the simple modules for all   Hecke algebras of type $G(d,p,n)$ in the modular case (see \cite{GJ}). Note also that results of  Hu provide  another parametrization of the simple modules for these algebras in the case where $p=d$ (see \cite{H} and \cite{H2}, see also \cite{J3} for the connections between these two classifications).

\bibliographystyle{amsalpha}

\end{document}